\documentclass[12pt,a4paper]{amsart}
\usepackage{amssymb}
\makeatletter
\renewcommand{\@begintheorem}[2]{
\rm \trivlist \item [\hskip \labelsep {\bf #2\ \ #1.}]
                                }
\makeatother

\makeatletter

\DeclareFontFamily{U}{cyr}{}
\DeclareFontShape{U}{cyr}{m}{n}{<5> wncyr5 <6> wncyr6 <7> wncyr7 <8> wncyr8 <9> wncyr9 <10->
wncyr10}{}
\DeclareMathAlphabet{\mathcyr}{U}{cyr}{m}{n}

\input cyracc.def

\newcommand{\ts}{\vspace{\baselineskip}\noindent{\bf Proof.}$\;\;$}
\newcommand{\ZZ}{{\bf Z}}
\newcommand{\QQ}{{\bf Q}}
\newcommand{\RR}{{\bf R}}

\newcommand{\PP}{{\bf P}}

\newcommand{\cC}{{\mathcal C}}
\newcommand{\cD}{{\mathcal D}}

\newcommand{\cF}{{\mathcal F}}

\newcommand{\cO}{{\mathcal O}}

\newcommand{\cS}{{\mathcal S}}
\newcommand{\cT}{{\mathcal T}}
\newcommand{\cU}{{\mathcal U}}
\newcommand{\cV}{{\mathcal V}}

\newcommand{\Pic}{\operatorname{Pic}}
\newcommand{\Pos}{\operatorname{Pos}}
\newcommand{\Mov}{\operatorname{Mov}}
\newcommand{\Nef}{\operatorname{Nef}}
\newcommand{\divi}{\operatorname{div}}
\newcommand{\Kthree}{\operatorname{K3}}
\newcommand{\bes}{\begin{equation*}}
\newcommand{\ees}{\end{equation*}}

\begin{document}




\date{}
\title[]
{The two extremal rays of some Hyper-K\"ahler fourfolds}

\author{Federica Galluzzi}
\address{Dipartimento di Matematica, Universit\`a di Torino, Via Carlo Alberto 10, Torino, Italy}
\email{federica.galluzzi@unito.it}

\author{Bert van Geemen}
\address{Dipartimento di Matematica, Universit\`a di Milano, Via Cesare Saldini 50, Milano, Italy}
\email{lambertus.vangeemen@unimi.it}


\begin{abstract}
We consider projective Hyper-K\"ahler manifolds of dimension four that are deformation equivalent to Hilbert squares of K3 surfaces. In case such a manifold admits a divisorial contraction, the exceptional divisor
is a conic bundle over a K3 surface. A classification of lattice embeddings implies that there are five types of such conic bundles.
In case the manifold has Picard rank two and has two (birational) divisorial contractions
we determine the types of these conic bundles. There are exactly seven cases. For the Fano varieties of cubic fourfolds there are only four cases and we provide examples of these.

\end{abstract}

\maketitle

\section*{Introduction}
The divisors on a four dimensional projective Hyper-K\"ahler manifold $X$ that can be contracted provide
interesting information on $X$ itself. In case $X$ is of K3$^{[2]}$-type and the Picard rank of $X$ is two,
as we will assume throughout this paper,
such a divisor $E$ is a $\PP^1$-fibration, so it is a conic bundle,
over a K3 surface $S$ of Picard rank one.
The fibration on $E$ is induced by the contraction of this divisor. The fourfold $X$ is then a moduli space
of possibly twisted sheaves on $S$ (\cite[Thm.\ 3.3]{vGeemenK}).
This implies that the Hodge structures on $H^2(X,\ZZ)$ and $H^2(S,\ZZ)$ are closely related.

The conic bundles on a K3 surface of Picard rank one which are isomorphic to such exceptional divisors were
classified in \cite{vGeemenK} and we recall this in Section \ref{sec:contractions}.
A main ingredient is the classification, due to Debarre and Macr\`i,
of embeddings of the rank two lattice $\langle 2d\rangle\oplus \langle -2\rangle$
into the lattice defined by $H^2(X,\ZZ)$ and its BBF-form.
An exceptional divisor defines an isometry class of such embeddings and there are five classes.

A well-known case is $E=\PP\cT_S$, the projectivization of the tangent bundle of a K3 surface $S$.
Then $X\cong S^{[2]}$, the Hilbert square of $S$. This type of contraction is denoted by H in this paper.
There are two other cases where $E=\PP \cV$ for a rank two vector bundle $\cV$ on $S$.
Then $\cV$ is a rigid, stable, vector bundle, called a Mukai bundle.
Such rank two bundles exist only if $h^2=2g-2$ with $g$ even, where $\Pic(S)=\ZZ h$.
We denote these cases by M$_i$ with $i=1,3$ according to $h^2\equiv 2i\mod 8$.
Even if the lattice embeddings of type $M_1$ and $M_3$ are quite different, the only notable difference
in the exceptional divisors is the genus of the K3 surface that is the base of the conic bundle.
Finally there are two cases where the Brauer class of the conic bundle is non-trivial.
In that case  $E=\PP\cU$ for a rigid $\alpha$-twisted locally free sheaf $\cU$ of rank two
and the Brauer class $\alpha\in Br(S)_2$ has order two.
These cases are denoted by B$_i$ with $i=0,1$ according to the intersection number
$Bh\equiv i/2\mod \ZZ$ where $B\in H^2(S,\QQ)$ is a B-field representative of $\alpha$.

An exceptional divisor on $X$ defines an extremal ray in $\Pic(X)_\RR\cong\RR^2$,
which is a boundary component of the movable cone of $X$.
The movable cone has two boundary components.
The other extremal ray is either generated by an isotropic vector, in which case it defines a Lagrangian fibration, or it defines another divisorial contraction of a HK manifold birational to $X$.
We recall the descriptions of the movable cone and the birational models of $X$ in Section \ref{sec:cones}.

In this paper we are interested in the types of these pairs of divisorial contractions.
A priori, given the 5 types of exceptional divisor, there could be $15$ types of such pairs.
However we show in Theorem \ref{thm:raypairs} that there are only seven cases.
Five of these are T+T, where T is any of the five types, the other two are H+M$_3$ and H+B$_1$.
Notice that H+M$_1$ cannot occur, illustrating the difference in types M$_1$ and M$_3$.
The types depend on rather subtle properties of the Picard lattice of $X$, like the divisibilities
in $H^2(X,\ZZ)$ of certain divisor classes associated to the contractions and the minimal solution of
a Pell equation associated to this lattice.
In the two `mixed' types H+M$_3$ and H+B$_1$, the presence of an H implies that $X$ is birational to
a Hilbert square $S^{[2]}$ for a K3 surface $S$.

The second exceptional divisor seems to play a rather important role in the geometry of the Hilbert square. The cases where $S$ has degree $14$ and $30$ the types are H+M$_3$ (see \S \ref{eis14} and \S \ref{Table2}).
The conic bundles $E$ of type M$_3$, which are scrolls over $S$, can be found in the literature
where also ad-hoc geometrical constructions of embeddings of $E$ in $S^{[2]}$ are given. However,
our result that such a conic bundle must be isomorphic to an exceptional divisor on $S^{[2]}$ may simplify and clarify the constructions.

In Section \ref{sec:Fano} we consider the special, but important, case that $X$ is birational to
the Fano variety of lines of a cubic fourfold. In that case $X$ can have only three types of contractions,
H, M$_3$ or B$_1$ and only four pairs of contractions are possible,
they are T+T, with T one of these three, or H+M$_3$.
We discuss the relation between the exceptional divisor and families of scrolls in the cubic fourfold,
already explored by Hassett and Tschinkel \cite{HassettT01} with an emphasis on the Mori cone of curves,
in terms of divisor classes.
We conclude with the well-known case of pfaffian cubic fourfolds.
Their Fano varieties are Hilbert squares $S^{[2]}$ and the degree of $S$ is $14$.
They have a contraction of type H and the second exceptional divisor is of type M$_3$.
As mentioned above, the corresponding conic bundle (it embeds as a scroll over $S$ in $\PP^5$)
was known and it was also shown to have an embedding in to $S^{[2]}$.
The identification of this conic bundle with (`the second') exceptional divisor seems to be new.






\

\section{Exceptional divisors on K3${}^{[2]}$ manifolds} \label{sec:contractions}

\subsection{Hyper-K\"ahler fourfolds of K3$^{[2]}$ type}\label{HK_K3^2}
A Hyper-Kähler (HK) variety is a simply-connected smooth projective variety admitting a
unique (up to scalar) non-degenerate holomorphic two-form.
These varieties are also called irreducible holomorphic symplectic varieties.
A HK fourfold $X$ is of K3$^{[2]}$ type if it is deformation equivalent to the Hilbert scheme of length 2 subschemes of a K3 surface.
The second cohomology group $H^2(X,\ZZ)$ is equipped with a quadratic form $q$, the Beauville-Bogomolov-Fujiki (BBF) form, we write $(x,y)$ for the associated bilinear form, so $q(x)=(x,x)$.
This quadratic form is related to the intersection product of divisors by the following relations,
$$
D^4\,=\,3q(D)^2\,=\,3(D,D)^2,\qquad D^3D'\,=\,3(D,D)(D,D')\qquad (D,D'\,\in\,\Pic(X))~.
$$

There is an isomorphism of lattices
$$
H^2(X,\ZZ)\,\cong\,\Lambda_{{\Kthree}^{[2]}}\,=\, \Lambda_{{\Kthree}}\oplus\,\ZZ\delta\,=\,
U^3\,\oplus\,E_8(-1)^2\,\oplus\,\langle -2\rangle~,
$$
where $\Lambda_{{\Kthree}}$ is the K3 lattice, $U$ is the hyperbolic plane and $E_8(-1)$ is the
unique unimodular, negative definite, even lattice of rank $8$.

The divisibility $\mbox{div}_{\Lambda_{\Kthree^{[2]}}}(v)$  of an element $v\in \Lambda_{\Kthree^{[2]}}$ is the positive integer such  that
$$
\{ (v,w):w\in \Lambda_{\Kthree^{[2]}}\}=\mbox{div}_{\Lambda_{\Kthree^{[2]}}}(v)\ZZ~.
$$
Any primitive $v\in \Lambda_{{\Kthree}}\subset \Lambda_{{\Kthree}^{[2]}}$ has divisibility one since
the K3 lattice is unimodular and hence is self-dual.
The divisibility of a primitive element in $H^2(F,\ZZ)\cong\Lambda_{\Kthree^{[2]}}$
is either one or two (\cite[Example 3.11]{Gritsenko_HS}).

\

\subsection{Contractions}\label{ssec:contractions}
Let $X$ be a complex manifold of K3${}^{[2]}$ type with Picard group of rank two.
If there exists a divisorial contraction $f:X\rightarrow Y$ then the exceptional
divisor $E$ is smooth and irreducible.
The class of $E$ in $\Pic(X)$ is either $\tau$ or $2\tau$ where $\tau\in \Pic(X)$
is a $(-2)$-class, that is, $q(\tau)=-2$.

The perpendicular of $E$ in $\Pic(X)$ is generated by a big and nef divisor class $H$.
A multiple of
$H$ defines the contraction, so one can identify $f$ with the map defined by $nH$ for some large $n$.
The contraction determined by $(H,E)$ defines a positive integer $d$ by $q(H)=2d$.
The Picard lattice then contains the sublattice
$$
\ZZ H\oplus \ZZ \tau\,\cong\,\langle 2d\rangle\,\oplus\,\langle -2\rangle\qquad(\subset\,\Pic(X))~.
$$
This sublattice is either equal to $\Pic(X)$ or it has index two in it, correspondingly one has
$|\det(\Pic(X))|=4d$ or $d$.

The classification of the contractions relies on the classification, due to Debarre and Macr\`\i$\;$
of the embeddings of lattices
$$
\ZZ H\oplus\ZZ \tau \;\hookrightarrow\;\Lambda_{{\Kthree}^{[2]}}~.
$$
In (\cite{Debarre_M}, \cite[Thm.\ 3.32]{Debarre}) it is shown that there are five classes of embeddings.

In these papers, the classes of the embeddings are related to
irreducible components of Heegner divisors in the boundary of moduli
spaces of polarized HK fourfolds of $\Kthree^{[2]}$-type
(in a general deformation of $(X,H)$, the big and nef divisor $H$ becomes ample).
The Heegner divisors corresponding to the five classes of embeddings are listed in
Table \ref{table:heegner} (cf.\ \cite[Table 1]{vGeemenK} where the order of the columns is different however)
but we will not consider these Heegner divisors in this paper.

For each of the five embeddings $\Pic(X)\cong \ZZ H\oplus\ZZ\tau\hookrightarrow H^2(X,\ZZ)$,
the isometry class of transcendental lattice $T(X):=\Pic(X)^\perp$
was determined in \cite{vGeemenK}.

\

\subsection{Conic bundles}\label{sec:conicbundles}
In \cite[Thm.\ 3.3]{vGeemenK} it is shown, with case by case computations, that a general
$X$ with a contraction defined by $(H,E)$ and Picard rank two is a moduli space of
(possibly twisted) sheaves on a K3 surface $S$.
The results of Bayer and Macr\`\i$\;$ in \cite{Bayer_M} then imply that
the contraction induces the structure of a conic bundle (a $\PP^1$-fibration) over $S$
on the exceptional divisor $E$:
$$
f:\,X\,\longrightarrow\,Y,\qquad p:=f_{|E}:\,E\longrightarrow\, S\;(\subset Y)~,
$$
and $f$ is an isomorphism on $X-E\cong Y-S$.
The Picard rank of $S$ is one and we write $\Pic(S)=\ZZ h_S$ where $h_S$ is ample. The self-intersection
number $h_S^2$ is either $q(H)$ or $q(H)/4$.

A conic bundle on $S$ defines a two-torsion Brauer class, denoted by $\alpha\in Br(S)_2$
where
$$
Br(S)\,:=\,H^2(S,\cO_S^\times)_{tors}\,\cong\, H^2(S,\QQ)/(\Pic(S)_\QQ\oplus H^2(S,\ZZ))~.
$$
A two-torsion Brauer class has a representative $B\in \mbox{$\frac{1}{2}$}H^2(S,\ZZ)$, called a B-field.
The intersection number $h_SB\in \mbox{$\frac{1}{2}$}\ZZ/\ZZ$ is an invariant of the Brauer class and,
only in the case that $4Bh_S+h_S^2\equiv 0\mod 4$, also $B^2\in \mbox{$\frac{1}{2}$}\ZZ/\ZZ$ is an invariant.
These invariants are also given in Table \ref{table:heegner}.
The transcendental lattice of $S$ is $T(S):=\Pic(S)^\perp$.
A Brauer class $\alpha\in Br(S)_2$ corresponds
to a homomorphism $\alpha:T(S)\rightarrow \mbox{$\frac{1}{2}$}\ZZ/\ZZ$ given by $t\mapsto Bt$.
The kernel of this homomorphism is denoted by $T(S,\alpha)$.
For non-trivial $\alpha$ this is an index two sublattice of $T(S)$.

\


\begin{table}
\caption{Divisorial contractions and Heegner divisors}
\label{table:heegner}
\centering
{\renewcommand{\arraystretch}{1.5}
\begin{tabular}{|l|c|c|c|c|c|c|}
\hline
Type&H&M$_1$&M$_3$&B$_1$&B$_0$\\
  \hline
Heegner div.\  &$\cD^{(1)}_{2d,2d,\alpha}$&$\cD^{(1)}_{2d,2d,\beta}$ &
  $\cD^{(2)}_{2d,2d,\alpha}$& $\cD^{(1)}_{2d,2d,\beta}$&$\cD^{(1)}_{2d,8d,\alpha}$  \\
\hline
$\mbox{div}_{H^2(X,\ZZ)}(H)$&1&1&2&1&1\\
\hline
$\mbox{div}_{H^2(X,\ZZ)}(\tau)$&2&1&1&1&1\\
\hline
$q(H)=2d,\quad d\equiv$&any&$1\mod 4$&$3\mod 4$&$0\mod 4$&any\\
\hline
$|\det(\Pic(X))|$&$4d$&$d$&$4d$&$4d$&$4d$\\
\hline
$T(X)$&$T(S)$&$T(S)$&$T(S)$&$T(S,\alpha)$&$T(S,\alpha)$\\
\hline
$h_S^2$&$2d$&$2d$&$2d$&$d/2$&$2d$\\
\hline
Invariants&0&0&0&$h_SB\,=\,\mbox{$\frac{1}{2}$},\;B^2=\mbox{$\frac{1}{2}$}$&
$h_SB\,=\,0,\;B^2=\mbox{$\frac{1}{2}$}$\\
  \hline
\end{tabular}
}
\end{table}

\subsection{The exceptional divisors and rigid rank two bundles}
The type of the embedding $\ZZ H\oplus\ZZ\tau\hookrightarrow \Lambda_{\operatorname{K3}^{[2]}}$
determines the invariants of the Brauer class of the conic bundle structure $p:E\rightarrow S$ on an exceptional divisor $E$, with class $E=\tau$ or $E=2\tau$, as well some other invariants.
In particular, there are also five types of conic bundles over a K3 surface that can appear as exceptional divisors.

The Brauer class of the conic bundle $E\rightarrow S$ is trivial for three types,
denoted by H, M$_1$ and M$_3$.
In these cases $E=\PP \cF$ where $\cF$ is a locally free sheaf of rank two on $S$.
In the case H, the contraction is the Hilbert Chow contraction
$S^{[2]}\rightarrow S^{(2)}$ and $\cF=\cT_S$, the tangent sheaf of $S$.
In the other two cases, $\cF=\cV$ where $\cV$ is the rank two rigid Mukai bundle on $S$.
If $\Pic(S)=\ZZ h_S$, such a bundle exists if and only if $h_S^2=2g-2$ with $g=2s$ even,
its Mukai vector is $v(\cV)=(2,h_S,2s)$ and rigidity is equivalent to $v(\cV)^2=-2$
(see \cite[\S 10.3.1]{Huybrechts_K3}). 
The $i$ in the type M$_i$ is determined by $h_S^2\equiv 2i\mod 8$.

If the Brauer class $\alpha\in Br(S)_2$ of $E\rightarrow S$ is non-trivial,
there are two types denoted by B$_0$ and B$_1$.
The exceptional divisor is then $E=\PP \cU$ for a rigid locally free rank two
$\alpha$-twisted sheaf $\cU$ on $S$ with Brauer class $\alpha$.
The $i$ in B$_i$ is determined by $h_SB\equiv i/2\mod \ZZ$
where $B$ is a B-field representative of $\alpha$.
In both cases one must have $B^2\equiv 1/2$, but recall that $B^2$ is not an invariant
unless $4Bh_S+h_S^2\equiv 0\mod 4$.
(On a general K3 surface $S$ there are three types of non-trivial Brauer classes in $Br(S)_2$.
In the case where $B^2$ is an invariant and $B^2\equiv 0$, there does not exist a locally
free rank two $\alpha$-twisted sheaf which is rigid 
and such Brauer classes are not obtained from exceptional divisors.)
Similar to the case $\alpha=0$, an $\alpha$-twisted sheaf $\cU$ has a (twisted) Mukai vector $v^B(\cU)$.
The rigidity of $\cU$ is equivalent to $v^B(\cU)^2=-2$ (cf.\ \cite[\S 2.3]{vGeemenK}).

Conversely, given a rigid, possibly twisted, locally free sheaf $\cU$ of rank two
on a K3 surface $S$ with Picard rank one as above, there is a HK fourfold $X$ of K3$^{[2]}$-type
that has an exceptional divisor $E$ with  conic bundle structure $E\cong\PP\cU\rightarrow S$.
This $X$ is a moduli space of twisted sheaves on $S$.

\subsection{Contractions of type H}\label{H}
For a contraction of type H the fourfold $X$ is isomorphic to a Hilbert scheme $S^{[2]}$ where $S$ is a K3 surface with $\Pic(S)=\ZZ h$, let $h^2=2d$.
The Picard lattice of $S^{[2]}$ is
$$
\Pic(S^{[2]})\,=\,\ZZ H\,\oplus\,\ZZ\tau,\qquad q(H)\,=\,2d,\quad (H,\tau)=0,\quad q(\tau)=-2,
$$
here $2\tau$ is the class of the exceptional divisor $E\cong \PP\cT_S$ in $S^{[2]}$ parametrizing the non-reduced subschemes of length two and $H$ is the divisor class induced by $h$.
The embedding $\ZZ H\oplus\ZZ\tau\hookrightarrow \Lambda_{\operatorname{K3}^{[2]}}$
is characterized, up to isometry,
by the fact that the divisibility of $\tau$ in $H62(X,\ZZ)\cong\Lambda_{{\Kthree}^{[2]}}$ is two.

\subsection{Contractions of type M$_1$}\label{M1}
According to Table \ref{table:heegner},
a contraction defined by $(X,H)$ is of type M$_1$ if and only if $q(H)=2d$ and $|\det(\Pic(X))|=d$.
Only in this case $\Pic(X)\neq\ZZ H\oplus\ZZ \tau$, but $\Pic(X)$ is generated by $H$ and $(H+\tau)/2$.
Since $q((H+\tau)/2)=(2d-2)/4$ must be an even integer, this implies that $d\equiv 1\mod 4$.
One can embed $\Pic(X)$ into  $\Lambda_{\mbox{\tiny K3}}\subset H^2(X,\ZZ)$, thus any
primitive element in $\Pic(X)$ has divisibility one in $H^2(X,\ZZ)$ (cf.\  \cite[Prop.\ 3.7]{vGeemenK}).


For example, if $d=1$, the even, indefinite, rank two Picard lattice is unimodular and hence is isometric
to the hyperbolic plane $U$.
The other extremal ray is then isotropic and thus defines a Lagrangian fibration,
these fourfolds are studied in \cite{DebarreHMV}.

\subsection{Contractions of type M$_3$}\label{M3}
A contraction defined by $(X,H)$ is of type M$_3$ if and only if the divisibility $\gamma$ of $H$ is two.
In that case, $H=2D+a\delta\in \Lambda_{\operatorname{K3}^{[2]}}$ for some $D\in\Lambda_{K3}$ and, since
$H$ is primitive, $a$ is odd. Then $q(H)=4D^2-2a^2\equiv -2\mod 8$, that is, $q(H)=2d$ with $d\equiv 3\mod 4$.
The Picard lattice is $\Pic(X)=\ZZ H\oplus\ZZ\tau$.
For examples, see \cite[Examples 5.4, 5.5, 5.6]{vGeemenK}.

\subsection{Contractions of type B$_0$} \label{B0}
The contractions with $\alpha\neq 0$, $Bh=0$ are the only ones where the transcendental lattice has a
non-cyclic discriminant group (\cite[Prop.\ 3.6]{vGeemenK}).
This implies that if $X$ has another
a divisorial contraction, then it has the same type. For these contractions, there is no congruence condition
on $d$.

In this case $\Pic(X)=\ZZ H\oplus\ZZ \tau$ and it can be embedded into $\Lambda_{\operatorname{ K3}}\subset H^2(X,\ZZ)$, see \cite[Prop.\ 3.6]{vGeemenK}.
Therefore any non-zero primitive element in $\Pic(X)$ has divisibility one in $H^2(X,\ZZ)$.

An example, due to O'Grady, is given in \cite[Prop.\ 5.1]{vGeemenK}.
In that example $d=1$ and the HK manifolds are double covers of EPW sextics,
the covering involution permutes the two exceptional divisors.

\subsection{Contractions of type B$_1$}\label{B1}
For a contraction with a non-trivial Brauer class $\alpha\in Br(S)_2$ and $Bh\equiv 1/2\mod\ZZ$,
the transcendental lattice $T(X)$ of $X$ is isomorphic to the index two sublattice $T(S,\alpha)$ of $T(S)$.
From  the description of the lattice  $T(S,\alpha)$ given in \cite[Prop.\ 3.5]{vGeemenK} one finds that its
discriminant group is cyclic and $q(H)=2d$ with $d\equiv 0\mod 4$.
For these $\alpha$ it can still happen that $T(S,\alpha)\cong T(S')$ for another K3 surface $S'$.
In fact, $T(X)\cong T(S')$, if and only if $d\equiv 0\mod 8$, see \cite[Prop.\ 3.5]{vGeemenK}.
In that case  $S'$ is a K3 surface of degree $d/2$.
Since $T(X)\cong T(S')\cong T({S'^{[2]}})$, the HK fourfolds $X$ and $(S')^{[2]}$ are called FM partners,
but in general this does not imply that $X$ and $(S')^{[2]}$ are birational.

An example, with $d=8$, is given in \cite[\S 5.3]{vGeemenK}, see also the case $e=2d=16$ in \S \ref{Table2}.

\

\subsection{Remark} From these descriptions, one finds that if $X$ has a contraction of type M$_1$ or B$_0$, then any other contraction of a birational model of $X$ will be of the same type.
The remaining three types, B$_1$, H, M$_3$, can be distinguished from each other by the divisibility
$\gamma$ of $H$ and the one of $\tau$:
\begin{itemize}
 \item[a)] type H if the divisibility of $\tau$ is two (and the divisibility of $H$ is one),
 \item[b)] type M$_3$ if the divisibility of $H$ is two (and the divisibility of $\tau$ is one),
 \item[c)] type B$_1$ if the divisibility of both $H$ and $\tau$ is one.
\end{itemize}

\

The following proposition, based on results of Druel, will be useful to identify exceptional divisors of
contractions.

\subsection{Proposition}\label{prop:excdiv}
Let $p:Z\rightarrow S$ be a conic bundle over a K3 surface and
assume that there is an embedding $Z\hookrightarrow X$ where $X$ is a HK fourfold.
Then there is a birational isomorphism $X\simeq X'$, with $X'$ a HK fourfold,
and the image of $Z$ is the exceptional fiber of a divisorial contraction $X'\rightarrow Y'$.

\ts
By \cite[Prop.\ 1.4, Theorem 1.3]{Druel} we only have to show that if $f\subset Z$ is a fiber of $p$ then
$Z\cdot f<0$. Let $i:Z\hookrightarrow X$ be the inclusion, then $Z\cdot f=Z\cdot i_*f_Z$ where $f_Z$ is
a fiber considered as a 1-cycle on $Z$. By the projection formula, this intersection number is $i^*Z\cdot f_Z=K_Z\cdot f_Z$
since $i^*Z=Z_{|Z}=K_Z$ by adjunction on $X$. As $K_Z=\omega_{Z/S}$ is also the relative dualizing sheaf
(cf.\ \cite[Section 1.1]{vGeemenK}) and $\omega_{Z/S}$ restricts to $\omega_f$ on $f$, we get $K_Z\cdot f_Z=\mbox{deg}(\omega_f)=-2<0$ as desired.
\qed

\

\section{The extremal rays of the movable cone} \label{sec:cones}

\subsection{Cones and divisorial contractions}
Recall that if $S$ is a K3 surface with Picard number two, then the divisorial contractions of $S$
are the contractions of a $(-2)$-curve $n\subset S$ and $n\cong\PP^1$ is a smooth rational curve.
The contraction $f:S\rightarrow \overline{S}$ of $n$ is defined by a big and nef divisor $h\in \Pic(S)$
with intersection number $hn=0$.
In this case a primitive embedding of the rank two Picard lattice of $S$ into the
(unimodular) K3 lattice is unique up to isometry.
In particular there are no issues with the divisibilities of the various divisor classes.
The half line $\RR_{\geq 0} h$ generated by $h$ in $\Pic(S)_\RR$ is a
boundary ray of the nef cone of $S$. Since $\Pic(S)_\RR\cong\RR^2$
and a cone in a two-dimensional real vector space has two extremal rays,
there is one other boundary ray.
This ray either defines another contraction,
in this case there are exactly two $(-2)$-curves in $S$, or it is $\RR_{\geq 0} e$ for an isotropic class
$e\in \Pic(S)$, so $e^2=0$, and then the ray defines a genus one fibration $S\rightarrow \PP^1$.
The latter case occurs if and only if the discriminant $|\det(\Pic(S))|$ of the Picard lattice is a square.

Let now again $X$ be a HK fourfold of K3$^{[2]}$-type with Picard rank two.
The positive cone $\Pos(X)$ of $X$ is
the connected component of the set $\{x\in \Pic(X)_\RR:\,q(x)>0\,\}$ that contains an ample class.
The nef cone $\Nef(X)$, which is the closure of the ample cone, is contained in the positive cone.
The movable cone $\Mov(X)$ is the (closed, convex) cone generated by classes of line bundles on $X$
whose base locus has codimension at least 2. There are inclusions:
$$
\Nef(X)\,\subset\,\Mov(X)\,\subset\,\Pos(X)\quad(\subset \Pic(X)_\RR\,:=\,\Pic(X)\otimes\RR)~.
$$
For a K3 surface, $\Nef(S)=\Mov(S)$.
Whereas a birational map between K3 surfaces is an isomorphism, this is not the case for HK fourfolds.
A birational map $\sigma:X\rightarrow X'$ between HK manifolds induces an isometry
$\sigma^*:H^2(X',\ZZ)\rightarrow H^2(X,\ZZ)$.
The movable cone $\Mov(X)$ is the union over all such birational maps $\sigma$
of the $\sigma^*(\Nef(X'))$.
These subcones are either equal or have disjoint interiors. The manifold $X$ has a unique smooth
birational model exactly when $\Mov(X)=\Nef(X)$.

We now recall how to determine these cones.

\

\subsection{The movable cone and the nef cone}
The interior of the movable cone is the connected component of
$$
\Pos(X)\,-\,\bigcup_{n\in \Pic(X),\atop q(n)=-2} \{x\in\Pos(X):\;q(x,n)=0\,\}
$$
that contains an ample divisor (\cite{Markman_Survey}, \cite[Thm.\ 3.16]{Debarre}). In particular,
if there are no $(-2)$-classes in $\Pic(X)$, then $\Mov(X)=\Pos(X)$.

If there are $(-2)$-classes,
the two boundary rays of $\Mov(X)$ are called the extremal rays of $X$.
Both extremal rays can be defined  by $(-2)$-classes
or one of these rays is defined by a $(-2)$-class and the other by a isotropic class in $\Pic(X)$
(so $q(x)=0$ for any $x$ on this ray).
The latter happens exactly when $|\det(\Pic(X))|$ is a square.
An extremal ray defined by $(-2)$-class defines a divisorial contraction (of birational models) of $X$
whereas an extremal ray spanned by an isotropic vector
defines a Lagrangian fibration on (a birational model of) $X$, see \cite{Markman}.

Recall that the interior of the nef cone is the ample cone.
If $\sigma:X\rightarrow X'$ is a birational map
between HK fourfolds then the interior of $\sigma^*\mbox{Nef}(X')$
is a connected component of
$$
\mbox{Int}(\Mov(X))\,-\,\bigcup_{\kappa\in \Pic(X),\, q(\kappa)=-10,
\atop
\divi_{H^2(X,\ZZ)}(\kappa)=2} \{x\in\Pos(X):\;q(x,\kappa)=0\,\}
$$
all connected components of this set are the ample cones of HK manifolds that are birational
to $X$ (\cite{Bayer_HT}, \cite[Thm.\ 3.16]{Debarre}).
The birational map $\sigma$ is a composition of Mukai flops
in (Lagrangian) planes. 

If there are no $(-10)$-classes with divisibility two, then $\Mov(X)=\Nef(X)$, otherwise at least one, and
possibly both, of the boundary rays of the nef cone are defined by such $(-10)$-classes and
these rays are called flopping walls.

Notice that if $\Pic(X)\subset \Lambda_{{\Kthree}}\subset H^2(X,\ZZ)$
then the unimodularity of
$\Lambda_{\mbox{\tiny K3}}$ implies that there cannot exist a (primitive) class of divisibility two in
$\Pic(X)$. In that case the movable cone is the nef cone.
This happens if $X$ has a contraction of type B$_0$ or M$_1$ (see \ref{B0}, \ref{M1}).


\

\subsection{The Pell equation}\label{pell}
Assume that $X$ has a divisorial contraction defined by $(H,E)$
and that $|\det(Pic(X))|=4d$ where $d$ is not a square.
Then $\Pic(X)=\ZZ H\oplus \ZZ\tau$ as before with $q(H)=2d$,
$q(\tau)=-2$ and $\RR_{\geq0}H$ is an extremal ray.
For a divisor class $yH+x\tau\in \Pic(X)$ one has
$$
q(yH+x\tau)\,=\,2dy^2-2x^2,\qquad\mbox{hence}\quad x^2\,-\,dy^2\,=\,1
$$
is the condition that $q(yH+x\tau)=-2$.
The Pell equation $x^2-dy^2=1$ has a positive solution, that is $(a,b)$ with integers $a,b>0$.
A positive solution with minimal $a$ is called minimal solution,
this is also the solution for which the slope $b/a$ is minimal (\cite[Appendix A]{Debarre}).
The following well-known lemma gives the second extremal ray in terms of the first one and the minimal
positive solution of the Pell equation.

\subsection{Proposition}\label{prop:second}
Let $(a,b)$ be the minimal positive solution of the Pell equation $x^2-dy^2=1$.
Then the other divisorial contraction is defined by $(H',\tau')$ with:
$$
H'\,:=\,aH\,-\,bd\tau,\qquad
\tau'\,:=\,bH\,-\,a\tau~.
$$

\ts
If $(a,b)$ is any solution of the Pell equation, then $\tau':=bH-a\tau$ is a $(-2)$-class in $\Pic(X)$ and
any $(-2)$-class is obtained from a solution. The divisor class $H':=aH-bd\tau$ is perpendicular to $\tau'$
and it is primitive since $a^2-db^2=1$ implies that $GCD(a,bd)=1$.

We claim that if a multiple of $H'$ is effective then $a,b>0$. In fact, since $H$ is big and nef, a positive multiple of $H^3$ is the class of a curve in $X$ which is not contained in $H'$
and thus $H^3H'\geq 0$.
Recall that $ H^3H'=3q(H)(H,H')$.
As $q(H)=2d>0$ one has $0\leq (H,H')=2ad$ hence $a\geq 0$. Similarly, $E^3=cf$ is a multiple of the
class of a fiber $f$ of the conic bundle $E\rightarrow S$ (\cite[Prop.\ 1.1]{vGeemenK}).
Since by adjunction $E_{|E}=K_E$ and
$K_E=\omega_{E/S}$, the relative canonical bundle, we get $Ef=K_Ef=-2$. Therefore $cEf=E^4=3q(E)^2>0$ hence $c<0$.
Since the class of $E$ is a positive multiple of
$\tau$, $f$ is a negative rational multiple of $\tau^3$.
Since $H'\neq E$ we have $fH'\geq 0$ and hence $0\geq \tau^3H'=3q(\tau)(\tau,H')=-6\cdot (2bd)$
so $b\geq0$ as we wanted.

Similarly, if a positive multiple of
$\tau'$ is effective then also $a,b>0$. Thus the second contraction is given by a pair $(H',E')$ with
$H',E'$ as in the lemma and $(a,b)$ is a positive solution of the Pell equation.

The description of the movable cone implies that the ray $\RR_{>0}H'$ is the one closest to $\RR_{>0} H$, that is with $b/a$ minimal, which concludes the proof of the lemma.
\qed

\

\section{Pairs of extremal rays} \label{sec:pairs}

\subsection{Theorem} \label{thm:raypairs}
Let $X$ be a HK fourfold of K3$^{[2]}$ type with Picard rank two, assume that $X$ admits a divisorial
contraction and that $|\det(Pic(X))|$ is not a square. Then the movable cone has two extremal rays
defining contractions
and $|\det(Pic(X))|\equiv 0,1\mod 4$.

The types of the extremal rays are:
\begin{enumerate}
\item[a)] In case $|\det(\Pic(X))|\equiv 1\mod 4$ the types are M$_1$+M$_1$.
\end{enumerate}
Now we assume that  $|\det(\Pic(X))|\equiv 0\mod 4$ and we write $|\det(\Pic(X))|=4d$.
\begin{enumerate}
\item[b)] If the discriminant group of $T(X)$ is not cyclic, the types are B$_0$+B$_0$.
\end{enumerate}
Now we assume that moreover the discriminant group of $T(X)$ is cyclic.
Let $(H,\tau)$, $(H',\tau')$ be the two pairs corresponding to the extremal rays.
Let $(x,y)=(a,b)$ be the minimal solution of the Pell equation $x^2-dy^2=1$.
\begin{enumerate}
\item[c)] In case 
$d\equiv 0\mod 4$ the types are:

\noindent
H+H if  div$_{H^2(X,\ZZ)}(\tau)=$div$_{H^2(X,\ZZ)}(\tau')=2$,

\noindent
B$_1$+B$_1$ if div$_{H^2(X,\ZZ)}(\tau)=$div$_{H^2(X,\ZZ)}(\tau')=1$,

\noindent
H+B$_1$ 
if 
$d\equiv 0\mod 8$ and $b$ is odd.
\item[d)] In case 
$d\equiv 1,2\mod 4$ the types are H+H.
\item[e)] In case 
$d\equiv 3\mod 4$ the types are

\noindent
H+H if div$_{H^2(X,\ZZ)}(\tau)=$div$_{H^2(X,\ZZ)}(\tau')=2$,

\noindent
M$_3$+M$_3$ if div$_{H^2(X,\ZZ)}(\tau)=$div$_{H^2(X,\ZZ)}(\tau')=1$,

\noindent
H+M$_3$ if $a$ is even.
\end{enumerate}

\ts
The Picard lattice has rank two and is even, therefore $|\det(Pic(X))|\equiv 0,1\mod 4$.

\noindent
a) In Remark \ref{M1} we observed that if $|\det(\Pic(X))|\equiv 1\mod 4$ then both rays have type M$_1$.

\noindent
b)
From Remark \ref{B0} we see that if  the discriminant group of $T(X)$ is not cyclic, the types are B$_0$+B$_0$.

If $|\det(\Pic(X))|\equiv 0\mod 4$ and the discriminant group is not cyclic,
there cannot be any rays of type M$_1$ or B$_0$. So we only need to consider the
remaining three types: H, M$_3$, B$_1$.

\noindent
d)
If there is a ray of type  B$_1$ then $|\det(\Pic(X)|=4d$ with $d\equiv 0\mod 4$ and if the type is M$_3$
then $d\equiv 3\mod 4$. So for $d\equiv 1,2\mod 4$ the types must be H+H.

\noindent
c)
For $d\equiv 0\mod 4$ one has types H+B$_1$ exactly when div$_{H^2(X,\ZZ)}(\tau)=2$ and
div$_{H^2(X,\ZZ)}(\tau')=1$ or vice versa.
To see that this corresponds to $b$ odd we use the embeddings of $\Pic(X)$ into $\Lambda_{\operatorname{K3}^{[2]}}=
U^3\oplus E_8(-1)^2\oplus \ZZ\delta$, $q(\delta)=-2$:
$$
\mbox{H}:\qquad H\,=\,(1,d)_1,\quad \tau\,=\,\delta
$$
(\cite[Prop.\ 3.4]{vGeemenK}) where $(r,s)_i$ lies in the $i$-th copy of $U$ and
$$
\mbox{B}_1:\qquad H\,=\,(1,d)_1,\quad \tau\,=\,(1,-d)_1+2(1,d/4)_2+\delta
$$
(\cite[Prop.\ 3.5]{vGeemenK}).
If the first extremal ray is of type H, then, by Proposition \ref{prop:second}, the second one has
$\tau'=bH-a\tau=b(1,d)_1-a\delta$. Thus $\tau'$ has divisibility one iff $b$ is odd.
If the first ray is of type B$_1$, then the second  one has
$\tau'=bH-a\tau=(b-a,(b+a)d)_1-2a(1,d/4)_2-a\delta$. The type of the second ray is H exactly when $\tau'$ has divisibility two, that is when $b-a$ is even.
Since $a,b$ cannot both the even because $a^2-db^2=1$, it follows that both must be odd.
Since $d\equiv 0 \mod 4$, $a$ is odd, hence the  type of the second ray is H if and only if $b$ is odd.

Finally we observe that if $a,b$ are odd, then $a^2-db^2\equiv 1\mod 8$ implies that $d\equiv 0\mod 8$.

\noindent
e)
Similarly, for $d\equiv 3\mod 4$ the types are either H+H, M$_3$+M$_3$ or H+M$_3$ and the latter
occurs if and only if the divisibilities of $\tau,\tau'$ are distinct. From \cite[Prop.\ 3.8]{vGeemenK}
the embedding can be chosen as
$$
\mbox{M}_3:\qquad H\,=\,2(1,(d+1)/4)_1+\delta,\qquad \tau\,=\,(1,-1)_2~.
$$
If the first extremal ray is of type H, then the second one has
$\tau'=bH-a\tau=b(1,d)_1-a\delta$ and $\tau'$ has divisibility $1$ if and only if $b$ is odd.
If the first ray is of type M$_3$, then the second one has
$\tau'=bH-a\tau=2b(1,(d+1)/4)_1+b\delta-a(1,-1)_2$ which has divisibility two if and only if $a$ is even.
Since $a^2-db^2=1$ and $d$ is odd, by going mod $2$, we see that $a$ even if and only if $b$ is odd.
\qed

\subsection{An application to Hilbert squares} \label{appl}
We consider the case $X=S^{[2]}$, the Hilbert square
of a K3 surface $S$ with $\Pic(S)\cong\ZZ h_S$ and $h_S^2=e$.
We denote by $H\in\Pic(S^{[2]})$ the divisor class corresponding to $h_S$, which has $q(H)=e$
(see also \cite[\S 3.7.2]{Debarre} for Hilbert schemes and their extremal rays).
The class $H$ is big and nef, it contracts $E=\PP \cT_S\subset S^{[2]}$. This contraction has type H.
The class of $E$ is
denoted by $2\tau\in \Pic(S^{[2]})$, one has $q(\tau)=-2$ and $q(H,\tau)=0$.
The classes $H,\tau$ are a basis of $\Pic(S^{[2]})$ and $|\det(\Pic(S^{[2]})|=2e$, so $e=2d$.
The divisibility of $H,\tau$ in $H^2(S^{[2]},\ZZ)$ is $1,2$ respectively.

Since one of the contractions has type H,  the cases (a), (b),
B$_1$+B$_1$ in (c) and
M$_3$+M$_3$ in (e) of Theorem \ref{thm:raypairs} do not occur.
The other type is thus determined by $e=2d$ and the parity of $b,a$.
(For Hilbert squares one can also use that the divisor class $rH+s\tau$ has divisibility $2$ in $H^2(X,\ZZ)$ exactly when $r$ is even and that $\tau'=bH-a\tau$ where $(a,b)$ is the minimal solution of $x^2-dy^2=1$.)
In Table \ref{table:Hilb2} we list the types and other useful information.

\

\begin{table}
\caption{Divisorial contractions on a Hilbert square of $S^{[2]}$, $S$ of degree $e$}
\label{table:Hilb2}
\centering
{\renewcommand{\arraystretch}{1.5}
\begin{tabular}{|l|c|r|}
$\begin{array}{ccccc}
\hline
e&(a,b)& \mbox{types} &H'&\tau'\\
\hline
2& & \mbox{H}&  &\\
\hline
4&(3,2) &\mbox{H+H}& 3H-4\tau& 2H-3\tau\\
\hline
6&(2,1)& \mbox{H+M}_3 & 2H-3\tau& H-2\tau\\
\hline
8&  &\mbox{H}&  & \\
\hline
10&(9,4) &\mbox{H+H}& 9H-20\tau& 4H-9\tau\\
\hline
12&(5,2) &\mbox{H+H}& 5H-12\tau& 2H-5\tau\\
\hline
14&(8,3) &\mbox{H+M}_3& 8H-21\tau& 3H-8\tau\\
\hline
16&(3,1) &\mbox{H+B}_1 & 3H-8\tau& H-3\tau\\
\hline
\end{array}$ &\qquad&
$\begin{array}{ccccc}
\hline
e&(a,b)& \mbox{types} &H'&\tau'\\
\hline
18&& \mbox{H}&  & \\
\hline
20&(19,6) & \mbox{H+H}& 19H-60\tau &6H-19 \tau\\
\hline
22& (10,3) & \mbox{H+M}_3& 10H-33\tau &3H-10\tau\\
\hline
24&(7,2) & \mbox{H+H}& 7H-24\tau &2H-7\tau\\
\hline
26&(649,180) & \mbox{H+H}& 649H-2340\tau & 180H-649\tau\\
\hline
28&(15,4) & \mbox{H+H}& 15H-56\tau & 4H-15\tau\\
\hline
30&(4,1) & \mbox{H+M}_3& 4H-15\tau & H-4\tau \\
\hline
32&$ $ & \mbox{H}& & \\
\hline
\end{array}$\\
\end{tabular}
}
\end{table}

\subsection{The Table \ref{table:Hilb2}}\label{Table2}
In Table \ref{table:Hilb2} we consider a K3 surface $S$ as in \S \ref{appl} and we determine the types of
the extremal rays of $S^{[2]}$.
In case $2e$ is a square, the second extremal ray defines a Lagrangian fibration.
Otherwise, the second extremal ray is spanned by a big and nef divisor $H'\in \Pic(S^{[2]})$.
The corresponding exceptional divisor is $E'$, with  class $2\tau'$ if this contraction
is of type H and it is $\tau'$ in the other cases. We give the minimal solution $(a,b)$ of the Pell equation
$x^2-dy^2$, $e=2d$, and we write $H',\tau'$ as linear combinations of $H,\tau$ following Proposition \ref{prop:second}.

Below we discuss some cases with small $e$.

\

In case $e=4$, $S\subset\PP^3$ is a quartic surface and $S^{[2]}$ has an involution (found by Beauville).
It is given by sending $p+q\in S^{[2]}$ to the residual intersection of the line spanned by $p,q$ with $S$.
This involution does not fix the exceptional divisor of non-reduced subschemes in $S^{[2]}$
and thus $S^{[2]}$ has two contractions of type H.

\

In the case $e=6$, the minimal solution of $x^2-3y^2=1$ is $(a,b)=(2,1)$
so the types are H+M$_3$, this example is discussed in \cite[\S 5.4]{vGeemenK}.
The equation $x^2-3y^2=5$ has no solutions modulo $3$, thus there are no flopping walls
and $S^{[2]}$ has a unique birational HK model.

By Proposition \ref{prop:second}, the other ray is spanned by $H':=2H-3\tau$
and the exceptional divisor $E'$, with class $\tau'=H-2\tau$, is a conic bundle $E'\cong\PP\cV$ with trivial
Brauer class where $\cV$ is the rank two Mukai bundle on $S$.

\

See Example \ref{eis14} for the case $e=14$.

\

In the case $e=16$, the minimal solution of the Pell equation $x^2-8y^2=1$ is $(a,b)=(3,1)$,
thus $b=1$ is odd and the types are H+B$_1$.
The equation $x^2-8y^2=5$ has no solutions modulo $8$ and thus there are no flopping walls.
The second exceptional divisor
$E'$ is a conic bundle over a K3 surface $S'$ of degree four.
It has non-trivial Brauer class $\alpha\in Br(S')_2$ and is known as the BOSS bundle,
cf.\ \cite[Prop.\ 5.3]{vGeemenK}, \cite{BOSS}.
There are Hodge isometries of transcendental lattices $T(X)\cong T(S)\cong T(S',\alpha)$.

\

In the case ${e=22}$, also discussed in \cite[\S 5.5]{vGeemenK},
the minimal solution of the Pell equation $x^2-11y^2=1$ is $(a,b)=(10,3)$ so the Hilbert scheme $X=S^{[2]}$
also admits (birationally) a contraction of type M$_3$ defined by the big and nef divisor $H$ with
$H'=10H-33\tau$. One has $q(H')=22$ and the divisibility of $H'$ is $2$.
The (singular) image of $S^{[2]}$ under map defined by $H'$ is a Debarre-Voisin variety,
see \cite[\S 3]{Debarre_V}.
In this case $X$ admits a flop, in fact $x^2-11y^2=5$ has the solution
$(a,b)=(7,2)$. The corresponding $(-10)$-class is $\kappa=2H-7\tau$ which has divisibility two. The orthogonal
of $\kappa$ is spanned by $7H-22\tau$ which lies in the movable cone, so the half line generated by the class
is a flopping wall.

\

The case $e=30$ is related to the recent preprint \cite[Conjecture 4.4.1(iv)]{Han} were, assuming a conjecture in the preprint, it is shown that $X:=\PP \cV$ embeds into $S^{[2]}$ where $\cV$ is the
rank two Mukai bundle on a general K3 surface $S$ of degree $e=30=2g-2$, so of genus $g=16$. From Proposition \ref{prop:excdiv} it follows
that the image of $X$ is an exceptional divisor.
There are no $(-10)$-classes in the Picard lattice, so $S^{[2]}$ has a unique HK model.
The results of this paper show directly that $X$ embeds in $S^{[2]}$. It is interesting to observe that $X$ appears to play a crucial role in showing that $S^{[2]}$ is DV variety.

\

\

\subsection{The K3 surfaces related to contractions}\label{FMpartners}
A divisorial contraction of a HK fourfold $X$ of $\Kthree^{[2]}$-type determines a K3 surface, which is the
base of the conic bundle structure $E\rightarrow S$ induced by the contraction,
where $E$ is the exceptional divisor.
From the classification of these contractions we see that  $T(X)\cong T(S)$ or $T(X)\cong T(S,\alpha)$.
However, these isomorphisms are in general not sufficient to determine $S$ uniquely.
In fact, two K3 surfaces with Hodge isometric transcendental lattices are called Fourier Mukai (FM) partners.
The number of FM partners of a given K3 surface is finite.
For example in \cite[Prop.\ 1.10]{Oguiso02} one finds that the number of FM partners
of a K3 surface $S$ with $\Pic(S)=\ZZ h$ and $h^2=2d$, is $2^{p(d)-1}$ where $p(d)$
the number of distinct primes dividing $d$.

Assume that $X$ has (birationally) two divisorial contractions, these define K3 surfaces $S_1,S_2$
with Picard rank one and such that $T(X)\cong T(S_i)$ for $i=1,2$. From Table \ref{table:heegner}
we see that the types must be H, M$_1$ or M$_3$. From Theorem \ref{thm:raypairs} we see that if the
types are not the same they must be H+M$_3$.
Let $\Pic(S_1)=\ZZ h_1$, $\Pic(S_2)=\ZZ h_2$, then again from Table \ref{table:heegner} we deduce that
$h_1^2=h_2^2=2d$ for some $d>0$. If $d$ is the power of a prime number then Oguiso's result implies
that $S_1\cong S_2$ but otherwise these K3 surfaces might not be isomorphic.
It is remarkable that $X$ then determines a specific pair of K3 surfaces whereas the condition
$T(X)\cong T(S)$ determines a possibly much larger set of $2^{p(d)-1}$ distinct K3 surfaces.

\


\section{Pairs of extremal rays of Fano fourfolds}\label{sec:Fano}

\subsection{The Picard lattice of a Fano fourfold}\label{mapalphan}
In this section, $Y$ will be a cubic fourfold. We denote by $F=F(Y)$ its Fano fourfold, which parametrizes the lines on $Y$, and which is a HK fourfold of K3$^{[2]}$ type.
There is an isomorphism (of abelian groups only, not of lattices!), \cite[\S 3, Prop.\ 6.3]{Beauville_D},
\cite[\S 2.5]{Huybrechts_cubicH},
$$
\alpha:\,H^4(Y,\ZZ)\,\stackrel{\cong}{\longrightarrow}\,
H^2(F(Y),\ZZ),
$$
which induces an isomorphism between the group $N^2(Y)$ of codimension two algebraic cycles on $Y$ and
$\Pic(F(Y))$.
This isomorphism does induce a Hodge isometry between the transcendental lattices
$$
\alpha:\,T(Y)(-1)\,\stackrel{\cong}{\longrightarrow}\,T(F),\qquad
\big(T(Y)\,:=\,N^2(Y)^\perp,\quad T(F)\,:=\,\Pic(F)^\perp\big)~,
$$
where the $T(F)(-1)$ indicates that the intersection form restricted to $T(Y)$ must be multiplied by $-1$.

Hassett showed that if $N^2(Y)$ has rank two, then the determinant $e$ of the intersection form on $N^2(Y)$
determines this lattice uniquely.
This lattice is denoted by $K_e$, one has  $e>6$ and $e\equiv0,2\mod 6$, so $e=2d$ is even.
The (irreducible) Hassett divisor $\cC_e$ in the moduli space of cubic fourfolds parametrizes the cubic fourfolds $Y$ such that $h_3^2\in K_e$, where $h_3\in H^2(Y,\ZZ)$ is the hyperplane class,
and $K_e$ is primitively embedded in $N^2(Y)$. The following proposition is well-known.

\

\subsection{Proposition}\label{PicFano}
Let $Y\in \cC_e$ be a cubic fourfold and assume that the rank of the lattice of
codimension two algebraic cycles is two. Let $F=F(Y)$ be the Fano fourfold of lines on $Y$,
and let $g\in Pic(F)$ the class of the Pl\"ucker polarization. Then $|\det(\Pic(F))|\,=\,2e=4d$
and there is a class $\gamma\in\Pic(F)$ such that the lattice $(\Pic(F),q)$ is
{\renewcommand{\arraystretch}{1.2}
$$
\Pic(F)\,=\,\left\{ \begin{array}{rcr}
\left(\ZZ g\,\oplus\ZZ \gamma,\,\begin{pmatrix}
              6&0\\0&-e/3
             \end{pmatrix}\right),&\qquad& \mbox{if}\quad e\equiv 0\mod 6~,
 \\ &&\\
\left(\ZZ g\,\oplus\ZZ \gamma,\,\begin{pmatrix}
              6&2\\2&(2-e)/3
             \end{pmatrix}\right),&\qquad& \mbox{if}\quad e\equiv 2\mod 6~.
             \end{array}\right.
$$
}
The divisibility in $H^2(F,\ZZ)$ of a primitive divisor class $rg+s\gamma \in\Pic(F)$ with $r,s\in\ZZ$ is:
$$
{\rm div}_{H^2(F,\ZZ)}(rg+s\gamma)\,=\,\left\{\begin{array}{rcl} 1&\mbox{if}&s\;\mbox{is odd},\\
                                               2&\mbox{if}&s\;\mbox{is even}~.
                                              \end{array}\right.
$$


\ts
Let $h_3\in H^2(Y,\ZZ)$ be the class of a hyperplane section, then $\alpha$ maps $h_3^2\in H^4(Y,\ZZ)$
to $g$ and $q(g)=2(h_3^2)^2=6$ where $q$ is the BBF form on $H^2(F,\ZZ)$.
The map $\alpha$ induces an isomorphism between $H^4(Y,\ZZ)_0$,
the perpendicular of $h_3^2$ in $H^4(Y,\ZZ)$, and $H^2(F(Y),\ZZ)_0$, the perpendicular of $g$,
which changes the sign of the intersection product.
As any element in $H^4(Y,\ZZ)$ can be written as $ah_3^2+bt$ with $t\in H^4(Y,\ZZ)_0$ and $3a,3b\in\ZZ$
and $g\cdot \alpha(ah_3^2+bt)=6a$ it follows that ${\rm div}_{H^2(F)}(g)=2$.

From the classification of the rank two lattices $N^2(Y)$ in \cite{Hassett_speccub} one has
{\renewcommand{\arraystretch}{1.2}
$$
N^2(Y)\,=\,\left(\ZZ h_3^2\,\oplus\,\ZZ R,K_e\right),
\quad\mbox{with}\quad K_e\,=\,\left\{ \begin{array}{ccr}
\begin{pmatrix} 3&0\\0&e/3\end{pmatrix}&\quad& \mbox{if}\quad e\equiv 0\mod 6~,\\
&&\\
\begin{pmatrix}
              3&1\\1&(e+1)/3
             \end{pmatrix}&\quad& \mbox{if}\quad e\equiv 2\mod 6~.
                    \end{array}\right.
$$
}

We already know that $q(g)=2\cdot 3=6$, where $g=\alpha(h_3^2)$, let $\gamma:=\alpha(R)$.

In case $e\equiv 0\mod 6$,
$R\in (h_3^2)^\perp$ and thus $q(\gamma)$ is the opposite of the self-intersection of $R$:
$q(\gamma)\,=\,-R^2\,=\,-d/3$. 

In case $e\equiv 2\mod 6$, $H^4(Y,\ZZ)_0\cap N^2(Y)$ is generated by $c=h_3^2-3R$ and $c^2=3-6+9\cdot(e+1)/3=
3e$. Therefore $\alpha(c)^2=-c^2=-3e$. Since $R=(h_3^2-c)/3$ we get
$$
q(\gamma)\,=\,(1/9)(\alpha(h_3^2)^2+\alpha(c)^2)\,=\,(1/9)(6-3e)\,=\,(1/3)(2-e)~.
$$
Similarly,
$$
g\cdot\gamma\,=\,(1/3)\alpha(h_3^2)(\alpha(h_3^2)-\alpha(c))\,=\,(1/3)(6+0)\,=\,2~.
$$


Recall that the divisibility of a primitive element in $H^2(F,\ZZ)\cong\Lambda_{\Kthree^{[2]}}$
is either one or two. 
The divisibility of $g$ is two, so if $s$ is even, the divisibility of $rg+s\gamma$ is two.
The lattice
$H^4(Y,\ZZ)$ contains the primitive sublattice generated by $h_3^2$ and $R$, hence it has
a $\ZZ$-basis $e_1=h_3^2$, $e_2=R$, $e_3,\ldots,e_{23}$. Since $H^4(Y,\ZZ)$ is unimodular (by Poincar\'e duality) there is a dual basis $e_1^*,\ldots,e_{23}^*$ with $e_i\cdot e_j^*=\delta_{ij}$. In particular,
$e_2^*\in H^4(Y,\ZZ)_0$ and $R\cdot e_2^*=1$. Since
$\gamma=\alpha(e_2)$ it follows that
$(rg+s\gamma)\cdot\alpha(e_2^*)=-s$ hence if $s$ is odd, the divisibility of $rg+s\gamma$ is one.
\qed

\

\subsection{K3 surfaces associated to cubic fourfolds}\label{sec:assK3}
If for $Y\in\cC_e$ the orthogonal complement to $K_e$ in $H^4(Y,\ZZ)$ is Hodge isometric (up to sign)
to the primitive cohomology of a K3 surface $S$ (of degree $e$)
then $S$ is said to be a K3 surface associated to the cubic fourfold $Y$ and $e$ is called admissible.
This is equivalent to the transcendental lattices $T(Y)(-1)$
and $T(S)$ being Hodge isometric.
Hassett showed that if $N^2(Y)=K_e$ then $Y$ has an associated K3 surface if and only if
the following condition $(**)$ is satisfied
(\cite[Thm.\ 5.1.3]{Hassett_speccub}, cf.\ \cite[\S 1.3]{Huybrechts_K3cat}):
$$
(**)\qquad \mbox{$e$ is even and $e/2$  is not divisible by nine or by any prime $p\equiv 2\mod 3$.}
$$
The first admissible $e$ are:
$$
e\,=\,14,\; 26,\; 38,\; 42,\; 62,\; 74,\; 78,\; 86,\;98,\;114,\;122,\;134,\;146,\;158,\ldots~.
$$
More generally, Huybrechts (\cite[Thm.\ 1.3]{Huybrechts_K3cat} showed that $T(Y)\cong T(S,\alpha)$
for a K3 surface $S$ and a Brauer class $\alpha\in Br(S)$ if and only if
the following condition (**)$'$ is satisfied
$$
\begin{array}{ll}
(**)'\quad &\mbox{$e$ is even and in the prime factorization $e/2=\prod p_i^{n_i}$  one has}
\\ &
\mbox{$n_i\equiv 0\mod 2$ for all primes $p_i\equiv 2\mod 3$.}
\end{array}
$$
If (**)$'$ is satisfied then the order of $\alpha\in Br(S)$ satisfies $\mbox{ord}(\alpha)^2|e$
(\cite[Lemma 2.13]{Huybrechts_K3cat}).


\

The following theorem provides the basic results on the pairs of extremal rays of a Fano variety.
In combination with Theorem \ref{thm:raypairs} it allows one to determine the types explicitly,
see Table \ref{table:Fano}.
The theorem also includes some well-known results on the K3 surfaces associated to Fano varieties.
See \cite[\S 6.5]{Huybrechts_cubicH} and \cite{Addington16}, \cite{Addington_Thomas} for further results
on the relations with (twisted) derived categories.

\subsection{Theorem} \label{thm:contrFano}
Let $F$ be the Fano variety of a cubic fourfold $Y\in\cC_e$ with $e=2d>6$ and assume that $\Pic(F)$
has rank 2. 
Then the following holds.

\begin{enumerate}
\item[a)] If there is an exceptional divisor on $F$, then
there is an isometry $\Pic(F)\cong<2d>\oplus <-2>$.
\item[b)]
A divisorial contraction of $F$ has type H, M$_3$ or B$_1$. If $F$ has a pair of extremal rays, then this pair
is one of H+H, H+M$_3$, M$_3$+M$_3$, B$_1$+B$_1$.
\item[c)]  
If $F$ has an extremal ray of type H or M$_3$, then  $e$ satisfies $(**)$ and
there is a K3 surface associated to $Y$.
There is an extremal ray of type H if and only if there exist integers $a,n$ such that $e=2(n^2+n+1)/a^2$.

\item[d)] If $F$ admits a divisorial contraction of type B$_1$ then $e$ satisfies $(**)'$ but not $(**)$,
in fact $e/4$ is an integer and satisfies $(**)$.
There is a twisted K3 surface $(S,\alpha)$ associated to $Y$ with $\alpha\in Br(S)_2$, $\alpha\neq 0$.

\item[e)] In case $e\equiv 0\mod6$ and $F$ has an extremal ray, then there is another extremal ray of the
same type.
In case the rays both have type H, $F$ has two divisorial contractions
inducing conic bundles $E_i\rightarrow S_i$ for K3 surfaces $S_1$, $S_2$
and $F\cong S_1^{[2]}\cong S_2^{[2]}$, so $F$ is ambiguous.
\end{enumerate}

\ts
a), b)$\quad$ Since the class $g\in\Pic(F)$ has divisibility two in $\Pic(F)$,
it is not the case that $\Pic(F)$ embeds into $\Lambda_{\Kthree}$. Hence extremal rays of type M$_1$ and B$_0$
are excluded. For the three remaining type of contractions we have seen that $\Pic(F)\cong <2d>\oplus <-2>$.

\noindent
c) For the Fano variety $F=F(Y)$ one has  $T(Y)\cong T(F)$.
From Table \ref{table:heegner} we see that if the type is H or M$_3$
then $T(F)\cong T(S)$ for some K3 surface $S$.
Hence $S$ is associated to $Y$.
The Fano $F$ has a ray of type H if and only if $F$ is birationally isomorphic to $S^{[2]}$
for a K3 surface $S$, this surface is associated to $Y$.
The existence of $a,n$ then follows from \cite[Prop.\ 6.1.3]{Hassett_speccub},
see \cite[Thm.\ 2]{Addington16}, \cite[Cor.\ 5.22]{Huybrechts_cubicH} for the converse.

\noindent
d) In case $F$ has a contraction of type B$_1$ then $d\equiv 0\mod 4$ hence $e$ does not satisfy $(**)$.
Moreover, $T(F)\cong T(S,\alpha)$ for a K3 surface $S$
and an $\alpha\in Br(S)$ of order two. 
As $T(Y)\cong T(F)$, the twisted K3 surface $(S,\alpha)$ is associated to $Y$.

We show that if $F$ allows a contraction, then either $e$ or $e/4$ satisfies $(**)$.
Notice that if $e/4$ satisfies $(**)$, then $e$ satisfies $(**)'$.
We verify that then the contraction is of type B$_1$.
We consider the two cases $e\equiv 0,2\mod 6$ separately, even if the proofs are rather similar.
Recall that for a prime number $p>3$ the congruence $x^2\equiv -3\mod p$ has a solution if and only if $p\equiv 1\mod 3$.

In case $e\equiv 2\mod 6$, let $k$ be such that $e=2+6k=2(1+3k)$ so that
$q(xg+y\gamma)=6x^2+4xy+(2-e)y^2/3=2(3x^2+2xy-ky^2)$.
Since $F$ has a contraction, there is a $(-2)$-class in $\Pic(F)$. So there are $x,y\in\ZZ$ such that $q(xg+y\gamma)=-2$,
equvalently $3q(xg+y\gamma)=-6$. Let $z=3x+y$ then this implies that
$$
z^2\,-\,(3k+1)y^2\,=\,-3
$$
has an integer solution. Let $p>3$ be a prime divisor of $e$, so of $3k+1$, then $z^2\equiv -3\mod p$
and thus $p\equiv 1\mod 3$. Obviously $3$ does not divide $e$. If $2^a$ divides $e$ and $a>3$, then
$e/2=3k+1\equiv 0\mod 8$, but $z^2\equiv -3\mod 8$ has no solutions, so $a\leq 3$. In case $a=2$,
$e/2=3k+1\equiv 2\mod 4$ and $z^2-2y^2\equiv -3\mod 4$ implies that $y$ is even, but then we find again
that $z^2\equiv -3\mod 8$, a contradiction. So $e=2^a\prod p_i^{m_i}$ where $a=1,3$
and all $p_i\equiv 1\mod 3$ which satisfies $(**)$ if $a=1$ and else $a=3$ and then
$e$ satisfies $(**)'$ while $e/4$ satisfies $(**)$.
In case $a=3$, one has the congruence $z^2-4y^2\equiv -3\mod 8$ which implies that $y$ is odd.
It follows that the divisibility of $xg+y\gamma$ in $H^2(F,\ZZ)$ is equal to one, hence the type cannot be H.
As $e\equiv 0\mod 4$, the type cannot be M$_3$ either, so it must be B$_1$.

In case $e\equiv 0\mod 6$, $q(xg+y\gamma)=6x^2-(e/3)y^2$.
Since we can write $e=6k$, there exists a $(-2)$-class iff $6x^2-2ky^2=-2$, that is
$$
3x^2\,-\,ky^2\,=\,-1~.
$$
Notice that $k$ obviously cannot be a multiple of $3$.
If a prime $p>3$ divides $k$ then we get $3x^2\equiv -1\mod p$ and thus $p\equiv 1\mod 3$.
Thus $e/6=k=2^{a-1}e'$ where $e'\equiv 1\mod 3$ and $a\geq 1$. Then we find $-2^{a-1}y^2\equiv -1\mod 3$
which shows that $a$ is odd. If $a>3$ we get $3x^2\equiv -1\mod 8$ which is impossible.
Thus $e=2^a\cdot3\cdot\prod p_i^{m_i}$ where $a=1,3$ and all $p_i\equiv 1\mod 3$.
If $a=3$ we have $k=e/6\equiv 4\mod 8$ and $3x^2-4y^2\equiv -1\mod 8$ implies that $y$ is odd.
Hence, as for the case $e\equiv 2\mod 6$, we get desired results.

\noindent
e) In case $e\equiv 0\mod 6$, one has $q(xg+y\gamma)=6x^2-(e/3)y^2$ and $g$ lies in the ample cone of $F$, hence
if $H=xg+y\gamma$ generates an extremal ray one has $y\neq 0$ and $H'=xg-y\gamma$ is also an extremal ray.
The corresponding $(-2)$-classes $\tau,\tau'$ are then of the form $rg\pm s\gamma$ and thus have the same
divisibility in $H^2(F,\ZZ)$ hence the types are the same by Theorem \ref{thm:raypairs}.
If there is a ray of type H,
$e$ satisfies $(**)$ hence $e/6\equiv 1\mod 3$ because any prime divisor of $e/6$ is $1$ mod $3$.
Then there are no $(-10)$-classes since $3x^2-(e/6)y^2=-5$ has no solution modulo $3$ hence $F$ has a unique birational HK model.
Therefore a ray of type H has an exceptional divisor $E\subset F$  with $E\cong\PP\cT_S$ for a
K3 surface $S$ and $F$ is isomorphic to $S^{[2]}$.
See also \cite[Prop.\ 6.2.2]{Hassett_speccub}.
\qed

\

\begin{table}
\caption{Divisorial contractions on the Fano variety of $Y\in\cC_e$}
\label{table:Fano}
\centering
{\renewcommand{\arraystretch}{1.5}
\begin{tabular}{|c|c|c|c|c|c|c|c|}
\hline
$e$&adm &$(-2)$-classes&$(a,b)$& types &$H$&$\tau$&deg$(\Sigma_f)$\\
\hline
$8$&$(**)'$ &Yes&$\Box$ & B$_1$& $g+\gamma$ & $\gamma$&2
\\
\hline
$14$&$(**)$ &Yes&$(8,3)$ &$\begin{array}{c}\mbox{H}\\ \mbox{M}_3\end{array}$&
$\begin{array}{c}3g+5\gamma\\ 3g-2\gamma \end{array} $&
$\begin{array}{c} g+2\gamma\\ g-\gamma \end{array}$&
$\begin{array}{c} 5 \\ 4 \end{array}$
\\
\hline
$24$&$(**)'$ &Yes&$(7,2)$ &$\begin{array}{c}\mbox{B}_1\\ \mbox{B}_1\end{array}$&
$\begin{array}{c}4g+3\gamma\\ 4g-3\gamma \end{array} $&
$\begin{array}{c} g+\gamma\\ g-\gamma \end{array}$&
$\begin{array}{c} 6 \\ 6 \end{array}$
\\
\hline
$26$& $(**)$&Yes&$(649,180)$& $\begin{array}{c}\mbox{H}\\ \mbox{H}\end{array}$&
$\begin{array}{c} 11g-7\gamma\\119g+137\gamma \end{array} $&
$\begin{array}{c} 3g-2\gamma\\ 33g+38\gamma \end{array}$&
$\begin{array}{c} 7 \\ 137 \end{array}$
\\
\hline
$38$& $(**)$&Yes&$(170,39)$& $\begin{array}{c}\mbox{H}\\ \mbox{M}_3\end{array}$&
$\begin{array}{c} 109g-61\gamma\\5g+4\gamma \end{array} $&
$\begin{array}{c} 25g-14\gamma \\ g+\gamma \end{array}$&
$\begin{array}{c} 61 \\ 8 \end{array}$
\\
\hline
$42$& $(**)$&Yes&$(55,12)$& $\begin{array}{c}\mbox{H}\\ \mbox{H}\end{array}$&
$\begin{array}{c} 14g+9\gamma\\ 14g-9\gamma \end{array} $&
$\begin{array}{c}  3g+2\gamma\\ 3g-2\gamma \end{array}$&
$\begin{array}{c} 9 \\ 9 \end{array}$
\\
\hline
$56$& $(**)'$&Yes&$(127,24)$& $\begin{array}{c}\mbox{B}_1\\ \mbox{B}_1\end{array}$&
$\begin{array}{c} 11g-5\gamma \\ 53g+37\gamma \end{array} $&
$\begin{array}{c} 2g-\gamma\\ 10g+7\gamma \end{array}$&
$\begin{array}{c} 10 \\ 74 \end{array}$
\\
\hline
$74$& $(**)$&No&$(6,1)$& None&None&None&None
\\
\hline
$78$& $(**)$&Yes&$(25,4)$& $\begin{array}{c}\mbox{M}_3\\ \mbox{M}_3\end{array}$&
$\begin{array}{c} 13g+6\gamma \\ 13g-6\gamma \end{array}$&
$\begin{array}{c}  2g+\gamma \\  2g- \gamma \end{array}$&
$\begin{array}{c} 12 \\ 12 \end{array}$
\\
\hline
\end{tabular}
}
\end{table}

\subsection{The Table \ref{table:Fano}}
In Table \ref{table:Fano} we list the contractions for the classical case $e=8$ (the cubic fourfold has a plane),
the cases $e=14,26,38,42,74,78$ which satisfy $(**)$ as well as the cases $e=24, 56$ which satisfy $(**)'$.

We use the basis $g,\gamma$ of $\Pic(F)$ as in Proposition \ref{PicFano}.
If a divisor $D=rg+s\gamma$ is effective, then, since $g$ is very ample,
we must have $g^3D>0$.
Proceeding as in the proof of Proposition \ref{prop:second} this gives
$$
0\,<\,g^3(rg+s\gamma)\,=\,3r(g,g)^2\,+\,3s(g,g)(g,\gamma)~.
$$
As $(g,g)=6$ and $(g,\gamma)=0,2$ for $e\equiv 0,2\mod 6$ respectively,
we see that $rg+s\gamma$ is an effective divisor if $r>0$ or $3r+s>0$ respectively .

To construct the table, we found the two pairs $(r,s)$ with $q(rg+s\gamma)=-2$, such that the lines
perpendicular to the divisor classes $rg+s\gamma$ define the chamber which contains the ample class $g$.
In case both extremal rays correspond to divisorial contractions we checked that the results from Proposition \ref{prop:second} hold. In case $d$ is a square, where $e=2d$, we write a $\Box$ in the column for the
minimal positive solution $(a,b)$ of the Pell equation $x^2-dy^2=1$.

\

In case $e\equiv 0\mod 6$ the BBF-form $q$ on the basis $g,\gamma$ of $\Pic(F)$ is in diagonal form,
$q(xg+y\gamma)=6x^2-(e/3)y^2$. The ample class $g$ lies in the movable cone and thus the two boundary,
extremal, rays of the movable cone are spanned by big and nef divisor classes $ng\pm m\gamma$ with $n>0$
and some $m$, similarly the $(-2)$-classes defining these rays are $\tau=rg\pm s\gamma$ with
$r>0$ and some $s$.

\

Each fiber $f$, a smooth rational curve, of the conic bundle $E\rightarrow S$ defines a scroll $\Sigma_f\subset Y$ in the cubic fourfold. This scroll is the union of the lines in $Y$ that correspond to the points in
$f\subset F(Y)$.
The degree of this scroll is the intersection number $g\cdot f$ which is the degree of the rational curve $f$
in the Grassmannian of lines in $\PP^5$ (cf.\ \cite[\S 7.1]{HassettT01}). It can be computed as follows.
Since $H$ contracts $E$ to $S$, we have $Hf=0$.
Recall that $Ef=-2$ (see the proof of Proposition \ref{prop:second}) and
that $\Pic(F)=\ZZ H\oplus\ZZ\tau$ with $q(\tau)=-2$.
In case $E=\tau$, a generator of $\Pic(F)$, we have $\tau f=-2$, otherwise the type is H and $E=2\tau$
and $\tau f=-1$. Therefore
$$
\mbox{deg}(\Sigma_f)\,=\,g\cdot f\,=\,\left\{\begin{array}{rcl} -2s&\mbox{if}& E\,=\,\tau~,\\
                                 -s&\mbox{if}& E\,=\,2\tau~,
                                \end{array}\right.\qquad \mbox{where}\;g\,=\,rH\,+\,s\tau~.
$$

\subsection{Example: $e=14$}\label{eis14}
Table \ref{table:Fano} shows that if ${\bf e=14}$ then there is a contraction of type H.
Thus the Fano variety of a general cubic fourfold $Y$ in $\cC_{14}$ is birational to $S^{[2]}$
for a K3 surface $S$ of degree $14$.
The transcendental lattice of such an $S$ is isomorphic to the one of $F$.
Since a K3 surface of degree 14 is uniquely determined by its transcendental lattice,
$S$ is uniquely determined by $F$ (that is, $S$ has no FM-partners, cf.\ \cite[Prop.\ 1.10]{Oguiso02}).
There are no $(-10)$-classes in $\Pic(F)$, hence $F$ is actually isomorphic to $S^{[2]}$.
See \cite{Beauville_D}, \cite[\S 6.2]{Huybrechts_cubicH} for this result and for
the geometrical relations between $S$, $Y$ and $F$.

Table \ref{table:Fano} also shows that there is a second contraction of $F$ of type M$_3$.
We recall some results that identify the exceptional divisor of this contraction.
There is an embedding $S\hookrightarrow Gr(1,5)$, let $\cU\rightarrow S$ be the restriction of the
universal bundle $\cS$ on $Gr(1,5)$ (cf.\ \cite{Beauville_D}, \cite[\S 3]{Mukai_88}).
The threefold $\PP\cU$ is identified in \cite[Corollary 6.2.16]{Huybrechts_cubicH} with a threefold
$\PP\cU\cong Z\subset\PP^5$.
(This threefold also appears in \cite{Beltrametti_SS}.)
There is an embedding $Z\hookrightarrow F$ (\cite[Remark 6.2.18]{Huybrechts_cubicH}).
From Proposition \ref{prop:excdiv} one deduces that $Z$
is the exceptional divisor of a contraction of $F$, this contraction is the one of type M$_3$.

The scrolls $\Sigma_f\subset Y$, where $f$ is a fiber of $Z\rightarrow S$, are well-known, see for example
\cite[6.2.6]{Huybrechts_cubicH}, they indeed have degree four (whereas those associated to $\PP\cT_S\subset S^{[2]}$ have degree five). Fano used these scrolls to prove the rationality of the general cubic fourfold in
$\cC_{14}$, cf.\ \cite[Theorem 4.3]{Alzati_R}.

From the description of contractions of type M$_3$ it follows that
$\cU\rightarrow S$ is, up to a twist by a line bundle, the Mukai bundle on $S$.
This bundle has Mukai vector $v=(2,h,4)$, so $v^2=-2$.
The stable bundle $\cU_v$ with this Mukai vector has $\det(\cU_v)=h$, $c_2(\cU_v)=5$ and $h^0(\cU_v)=6$.
As $H^0(\PP{\cU_v},\cO_{\PP\cU_v}(1))\cong H^0(S,\cU_v)$ we obtain a map from $\PP{\cU_v}$ to $\PP^5$,
its image should be $Z$.
The kernel of the  evaluation map $H^0(\cU_v)\otimes\cO_S\rightarrow\cU_v$
will give the `Gushel' embedding $S\hookrightarrow Gr(4,H^0(\cU_v))\cong Gr(1,5)$ (see also \cite{Bayer_KM}).

\

\subsection{Rationality of cubic fourfolds} The cubic fourfolds in the Hassett divisors $\cC_e$, with
$e=14,26,38,42$ are known to be rational (established by Morin and Fano for $e=14$ and by Russo and Staglian\'o \cite{Russo_St19}, \cite{Russo_St23} in the other cases).
The types of the extremal rays of the corresponding Fanos do not seem to have an obvious pattern and it is not clear that the knowledge of the types of the extremal rays is useful for establishing the rationality. However, the proofs of the rationality often use scrolls in the
cubic fourfolds obtained from fibers of the conic bundle structure $p:E\rightarrow S$ of an exceptional divisor.
Notice that the cubic fourfolds $Y$ in $\cC_e$ with $e=74$ are conjectured to be rational,
but the Fano variety of the general $Y$ in $\cC_{74}$ does not have any extremal ray.

\

\end{document}